\documentclass[12pt,a4paper]{article}

\usepackage{amssymb,amsmath,amsthm,xspace,ScriptFonts,epsfig}

\begin{document}

\newcommand{\RR}{\mathsf{R}}
\newcommand{\LL}{\mathsf{L}}

\newcommand{\R}{\Bbb R}
\newcommand{\Z}{\Bbb Z}
\newcommand{\Q}{\Bbb Q}
\newcommand{\C}{\Bbb C}
\newcommand{\esp}{\vskip .3cm \noindent}
\newcommand{\el}{{\cal L}}
\newcommand{\de}{{\cal D}}
\newcommand\Cfty{${\cal C}^{\infty}$}
\newcommand{\emm}{{\cal M}}
\newcommand{\be}{{\cal B}}
\newcommand{\ess}{{\cal S}}
\newcommand{\enn}{{\cal N}}
\newcommand{\ka}{{\cal K}}
\newcommand{\emmh}{\hat {\cal M}}
\newcommand{\RD}{\text{\rm RD}}
\mathchardef\flat="115B
\newcommand{\trf}{\text{Tr}^\flat}
\newcommand{\dets}{\text{Det}^\#}
\newcommand{\detn}{\text{Det}^*}
\newcommand{\detf}{\text{Det}^\flat}
\newcommand{\trn}{\text{Tr}^*}
\newcommand{\trs}{\text{Tr}^\#}
\newcommand{\definition}{\overset{\text{\rm def}}{=}}

\newcommand\adm{{\subset}^a}
\newcommand\cfp{cofinal${}^+$}
\newcommand\cfm{cofinal${}^-$}
\newcommand\bdp{bounded${}^+$}
\newcommand\bdm{bounded${}^-$}

\def\p#1{ \buildrel \text{.} \over #1 }
\def\pp#1{ \buildrel \text{..} \over #1 }
\def\ppp#1{ \buildrel \text{...} \over #1 }
\def\ut#1{$\underline{\text{#1}}$}
\def\CC#1{${\cal C}^{#1}$}
\def\h#1{\hat #1}
\def\t#1{\tilde #1}
\def\wt#1{\widetilde{#1}}
\def\wh#1{\widehat{#1}}
\def\wb#1{\overline{#1}}

\def\restrict#1{\bigr|_{#1}}

\newtheorem{thm}{Theorem}

\newtheorem{thmalph}{Theorem}
\renewcommand\thethmalph{\Alph{thmalph}}

\newtheorem{ex}{Example}[section]

\newtheorem{prop}{Proposition}[section]
\newtheorem{lemma}[prop]{Lemma}
\newtheorem{cor}[prop]{Corollary}
\newtheorem{rem}[prop]{Remark}
\newtheorem{defi}[prop]{Definition}

\newtheorem*{prop*}{Proposition}
\newtheorem*{claim}{Claim}
\newtheorem*{lemma0}{Lemma 0}

\title{The homotopy classes of continuous maps\\
between some non-metrizable manifolds}
\author{Mathieu Baillif}
\date{\empty}

\maketitle

{\bf keywords:} {\em homotopy, longline, nonmetrizable manifolds}\\

{\bf Subject code classification:} {57N99, 14F35}

\abstract
We prove that the homotopy classes of continuous maps $\RR^n\to\RR$, where $\RR$ is Alexandroff's long ray, 
are in bijection with the antichains of $\mathcal{P}(\{1,\dots,n\})$. The proof uses partition
properties of continuous maps $\RR^n\to\RR$. We also provide a description of $[X,\RR]$ for some
other
non-metrizable manifolds.
\endabstract
\section{Introduction}
This paper is about computation of homotopy classes of maps between some non-metrizable manifolds.
The main result is a complete classification of homotopy classes of continuous functions
$\RR^n\to\RR$,
where $\RR$ is
Alexandroff's long ray,
which are shown to be in bijection with the antichains of
$\mathcal{P}(\{1,\dots,n\})\backslash\{\emptyset\}$
(see Theorems \ref{thm1}--\ref{thm2} below). This generalises a result
of D. Gauld \cite{Gauld} who solved the problem when $n=1$.

Under our opinion, there are (at least) three reasons that motivate this investigation of homotopy in non-metrizable manifolds. Firstly,
any manifold is compactly generated, and thus, according to G.W. Whitehead,
fits in the natural category of homotopy theory (see e.g. \cite{Whitehead}).
Secondly,
our manifolds provide a class of spaces $X$ for which $\Pi_i(X)=0$ for each $i\in\omega$ while
$[X,X]$ is finite but has at least two elements, so in particular $X$ is not contractible.
(Notice that since $\Pi_1(X)=\{0\}$, we do not need to bother about base points, and consider only free
homotopy.)
The non-contractibility does not come from the ``shape'' of $X$ but rather from its ``wideness''.
The third reason (of a more practical nature) is that the
proofs are completely elementary, in the sense that we use only very basic facts
about countable ordinals.
In fact, despite what the title may suggest, the main part of this paper consists in an investigation
of partition properties of maps $\RR^n\to\RR$ (which we find
interesting in themselves), which enable us to reduce the purely homotopical
questions to the trivial fact that two maps $[0,1]^n\to[0,1]$ are homotopic.

The paper is organised as follows. Section 2 contains the definitions and the statements of the main results. 
In particular, we define the cofinality class $\mathfrak{C}(f)$ of a map
$f:\RR^n\to\RR$. In Section 3, we devise some properties of ``big'' open and closed sets in $\RR^n$
which will often be of use. In Section 4, we show that the cofinality classes of functions
$\RR^n\to\RR$ are in bijection with the antichains of
$\mathcal{P}(\{1,\dots,n\})\backslash\{\emptyset\}$. Then, in Sections 5-6, we prove that $f$ and $g$
are homotopic if and only if $\mathfrak{C}(f)=\mathfrak{C}(g)$. Finally, in Section 7 we investigate some other
non-metrizable manifolds.
\esp
This paper can be seen as a companion to \cite{meszigues+Dave} where D. Cimasoni and I have
investigated
embeddings $\RR\to\RR^n$ up to ambiant isotopy.


\section{Definitions}
We recall that Alexandroff's (closed)
long ray $\RR$ is $\omega_1\times [0,1[$ endowed with the topology
given by the lexicographic order $\le$.
It is well known that $\RR$ can be made into a $1$-dimensional
(\CC{\infty}) manifold, is sequentially compact, non-metrizable and non-contractible.
In this paper, sequential compacity
is the key property and will always be implicitely invoked when we say that some (sub)sequence
converge. Two other well known properties of $\RR$ are given in the following lemmas whose proofs
can be found e.g. in \cite[Lemma 3.4 (iii)]{Nyikos:1984} and \cite{Kunen}:

\begin{lemma}\label{lemmaintro2}
  Let $f:\RR\to\RR$ be continuous and bounded. Then, $f$ is eventually
  constant, i.e. there is $z\in\RR$ such that $f(x)=f(z)$ if $x\ge z$.
\end{lemma}

\begin{lemma}\label{lemmaintro1}
  Let $\{E_m\}_{m\in\omega}$ be closed unbounded sets of $\RR$. Then,
  $\cap_{m\in\omega}E_m$ is closed and unbounded.
\end{lemma}

(In both lemmas $\RR$ can be replaced by $\omega_1$.)
We will always identify the ordinal $\alpha\in\omega_1$ with $(\alpha,0)\in\RR$, and thus consider $\omega_1$
as a subset of $\RR$. We use greek letters for ordinals, and only for them.\\
Let us fix $n$, the dimension, and set $N=\{1,\dots,n\}$. We denote by $\pi_i:\RR^n\to\RR$ the projection on the $i$-th
coordinate. We will often define sequences in $\RR$ or $\RR^n$, so to avoid confusion
we shall use only the index $m$ to denote a member of a sequence while we reserve the indices
$i,j,k,\ell$ for coordinates. For $x=(x_1,\dots,x_n)\in\RR^n$ we set $|x|=\max_{i=1,\dots,n}x_i$.
For a finite set $I$, we denote its number of
elements by $|I|$. If $I=\{i_1,\dots,i_k\}\subset N$ and $x\in\RR^n$, we write $x_I$ for
$(x_{i_1},\dots,x_{i_k})$.

\begin{defi}
  Let $I\subset N$ and $c\in\RR^{n-|I|}$.
  The $I$-diagonale at height $c$ is the set
  \begin{equation}
    \Delta_I(c)\definition\left\{x=(x_1,\dots,x_n)\in\RR^n\, : \,
    \begin{array}{c}
      x_i=x_{i'} \text{ if } i,i'\in I \\
      \text{ and } x_{N\backslash I}=c
    \end{array}
    \right\}.
    \label{I-diagonale}
  \end{equation}
  We abbreviate $\Delta_I(0)$ by $\Delta_I$.
\end{defi}
Notice that $\Delta_I(c)$ is homeomorphic to $\RR$.

\begin{lemma}\label{lemma0} Let $I\subset N$ and $c,c'\in\RR^{n-|I|}$.
  Then, $f\bigr|_{\Delta_I(c)}$ is unbounded (resp. bounded) if and only
  if $f\bigr|_{\Delta_I(c')}$ is unbounded (resp. bounded).
\end{lemma}

\proof
  Let $\gamma:[0,1]\to \RR^{n-|I|}$ be continuous with $\gamma(0)=c$, $\gamma(1)=c'$.
  Then, $\gamma$ provides an homotopy between $f\bigr|_{\Delta_I(c)}$
  and $f\bigr|_{\Delta_I(c')}$. Thus, since $\RR$ is non-contractible,
  both are either unbounded or bounded.
\endproof
One checks easily that if $I\not= J$, there is no homotopy sending $\Delta_I$ to $\Delta_J$; this and Lemma \ref{lemma0}
motivate the following definition:
\begin{defi}
  Let $f:\RR^n\to\RR$ be continuous and $I\subset N$.
  We say that $f$ is $I$-cofinal (resp. $I$-bounded)
  if $f\bigr|_{\Delta_I}$ is unbounded (resp. bounded).
\end{defi}

\begin{defi} The cofinality class of $f:\RR^n\to\RR$
  is the set $\mathfrak{C}(f)\definition\{I\subset N\,:\, f \text{ is } I\text{-cofinal}\}$.
\end{defi}
Our main results are :

\begin{thm}\label{thm1}
  Two continuous maps $f,g:\RR^n\to\RR$ are homotopic if and only if $\mathfrak{C}(f)=\mathfrak{C}(g)$.
\end{thm}

We recall that an antichain in a partially ordered set is a set of pairwise incomparable
elements. As usual, we order $\mathcal{P}(N)$ by the inclusion.
\begin{thm}\label{thm2}
  The homotopy classes $[\RR^n,\RR]$ of continuous maps
  $\RR^n\to\RR$ are in bijection
  with the antichains
  of $\mathcal{P}(N)\backslash\{\emptyset\}$.
\end{thm}

It is worth noting that the problem of counting the antichains of $\mathcal{P}(N)$ is
{\sc NP}-complete, see \cite{Ball+Provan}. The exact values for $n=1,\dots,7$ as well as some
inequalities can be found in \cite{Anderson}.


\section{Topology in $\RR^n$}
We first prove a useful property of ``big'' open sets, that is, those who contain some $\Delta_I(c)$
outside of a compact set. The formulation given below
is slightly more general.

\begin{lemma}[structure of open sets]
  \label{structure open sets}
  Let $h:\RR\to\RR^n$ be continuous such that $\pi_i\circ h$ is cofinal for
  $i\in I$ and bounded otherwise, and let $U$ be an open set
  in $\RR^n$ containing $\text{\rm Im}h$. Then, there are $x\in\RR$ and $y_j,y_j'\in\RR$
  ($j\in N\backslash I$) with
  $y_i<y_i'$, such that
  $$
    U\supset \prod_{i=1}^n U_i,
  $$
  where $U_i=[x,\omega_1[$ if $i\in I$ and $U_j=]\,y_j,y_j'[\,$ if $j\notin I$.
\end{lemma}
Notice that $x$ does not depend on $i$.
Given $z\in\RR$, it is easy to find $h$ such that $\text{\rm Im}h=\Delta_I(c)\backslash\,]0,z[^n$.
\proof
  We may assume that $I=\{1,\dots,s\}$.
  Since $\pi_i\circ h$ is bounded for $i=s+1,\dots,n$, there are $z,c_{s+1},\dots,c_n\in\RR$
  such that $\pi_i\circ h\bigr|_{[z,\omega_1[}\equiv c_i$.
  Write $c=(c_{s+1},\dots,c_n)$. Suppose (ab absurdo) that
  for all $x,y_i,y_i'\in\RR$ with $y_i<y_i'$,
  $$
    [x,\omega_1[^s\times\left(\prod_{k=s+1}^n]y_i,y_i'[\right)\not\subset U.
  $$
  We shall show that
  this implies that $\Delta_{\{1,\dots,s\}}(c)\cap(\RR^n\backslash U)$ is (closed and) unbounded
  (in $\Delta_{\{1,\dots,s\}}(c)$), which is a contradiction since
  $\text{\rm Im}h\cap\Delta_{\{1,\dots,s\}}(c)$ is (closed and) unbounded as well.
  (It is well known that if $g:\RR\to\RR$ is continuous and cofinal, the set of its fixed points
  $\{x:g(x)=x\}$ is closed unbounded. Apply Lemma \ref{lemmaintro1} to the
  $s$ first coordinates of $h$.)
  \\
  So, let $u\in\RR$. For $i=s+1,\dots,n$, let us fix sequences
  $y_{i,m}\nearrow c_i$ and $y_{i,m}\searrow c_i$ ($m\in\omega$), and set
  $x_0=(u,\dots,u)\in\RR^s$. We choose by induction the sequences
  $x_m\in\RR^s$ and $z_m=(z_{s+1,m},\dots,z_{n,m})\in\RR^{n-s}$ such that:
  \begin{equation*}
    x_m\in[\,|x_{m-1}|,\omega_1 \,[^s, \quad
    z_{i,m} \in\,\, ]y_{i,m},y_{i,m}'[,\quad\text{and } (x_m,z_m)\notin U
  \end{equation*}
  Then, $z_{i,m}\to c_i$, and $x_m\to(x,\dots,x)$ for some $x\ge u$. By closeness,
  $(x,\dots,x,c)\in\Delta_{\{1,\dots,s\}}(c)\cap\RR^n\backslash U$.
\endproof

We now prove a kind of analog of Lemma \ref{structure open sets} for closed sets which will
be useful in the proof of Lemma \ref{partition lemma}.
\begin{lemma}[structure of closed unbounded sets]\label{structure closed sets}
  Let $F\subset\RR^n$ be closed and unbounded. Then, there are $I\subset N$ and
  $c\in\RR^{n-|I|}$ such that $\Delta_I(c)\cap F$ is closed and unbounded
  (in $\Delta_I(c)$).
\end{lemma}
\proof
  Let $U=\RR^n\backslash F$ and $J=\{i\in N\,:\,\pi_i(F)\text{ is unbounded}\}$. We may
  assume $J=\{1,\dots,s\}$. For $i=s+1,\dots,n$, $\pi_i(F)$ is bounded by $b\in\RR$, say.
  Thus, $F\subset\RR^s\times[0,b]^{n-s}$. We show the result by induction on $s=|J|$.
  Suppose that for all $c=(c_{s+1},\dots,c_n)\in[0,b]^{n-s}$,
  $\Delta_{J}(c)\cap F$ is bounded. Thus, for each such $c$ there is
  $x(c)$ such that $U\supset (\Delta_J\backslash[0,x(c)[^n)$. By Lemma
  \ref{structure open sets}, there are $x'(c)$ and $y(c_i)<c_i<y'(c_i)$ such that
  $$
    [x'(c),\omega_1[\times\Bigl(
    \underbrace{\prod_{i=s+1}^n\,]y(c_i),y'(c_i)[\,}_{V(c)} \Bigr)\subset U.
  $$
  Since $\{V(c)\}_{c\in[0,b]^{n-1}}$ is an open cover of $[0,b]^{n-1}$, there is a finite
  cover $\{V(c^1),\dots,V(c^m)\}$. Putting $x=\max_{j=1,\dots,m}x'(c^j)$, we get
  $U\supset [x,\omega_1[^s\times [0,b]^{n-s}$. Letting
  $Q_i=\{(x_1,\dots,x_s)\in\RR^s\,:\, x_i\in[0,z]\}\times[0,b]^{n-s}$, we proved
  that $F\subset\cup_{i=1}^s Q_i$ and then for some $i\in\{1,\dots,s\}$,
  $F\cap Q_i$ is unbounded. But the number of unbounded projections of $F\cap Q_i$ is
  at most $s-1$ (which is a contradiction if $s=1$ since $F$ is unbounded),
  and we finish by induction.
\endproof


\section{Cofinality classes}

\begin{lemma}\label{lemma cof 1} Let $I\subset N$ be non-empty.
  If $f$ is $I$-cofinal and $J\supset I$, then $f$ is $J$-cofinal.
\end{lemma}

\proof 
  We may assume that $I=\{1,\dots,s\}$ and $J=\{1,\dots,s'\}$ with $s'\ge s$.
  Let $z\in\RR$, we shall find $x\in\Delta_J$ satisfying $f(x)\ge z$.  
  Put $x_0=0$ and define $x_m\in\RR$ ($m\in\omega$) as follows.
  Given $x_{m-1}$, take $x_m\ge x_{m-1}$ such that
  $$
    f(\underbrace{x_m,\dots,x_m}_{s} ,\underbrace{x_{m-1},\dots,x_{m-1}}_{s'-s},
    \underbrace{0,\dots,0}_{n-s'})\ge z
  $$
  ($x_m$ exists since by Lemma \ref{lemma0}, 
  $f\bigr|_{\Delta_I(x_{m-1},\dots,x_{m-1},0,\dots,0)}$
  is unbounded).
  This sequence converge to some $x$ and we have
  $$ 
    f( \underbrace{x,\dots,x}_{s'},\underbrace{0,\dots,0}_{n-s'})\ge z.
  $$  
\endproof

\begin{lemma}\label{lemma cof 2}
  If $f$ unbounded then $f$ is $N$-cofinal.
\end{lemma}

\proof
  By induction on $n$. If $n=1$, the result is trivial.
  So, let $n>1$ and suppose that $f\bigr|_{\Delta_N}$ is bounded, say by $b\in\RR$.
  Fix $b'>b$ and let
  $F$ be the closed unbounded set $f^{-1}([b',\omega_1[)$. Then, the open set $U=\RR^n\backslash F$
  contains $\Delta_N$.
  By Lemma \ref{structure open sets}, there is some $z\in\RR$ such that
  $U\supset [z,\omega_1[^n$.
  Thus, $F\subset\cup_{i=1}^nQ_i$, where
  $Q_i=\{(x_1,\dots,x_n)\in\RR^n\,:\, x_i\in [0,z]\}$,
  and $f\bigr|_{Q_i}$ is unbounded for some $i$. We fix this $i$.
  For $c\in\RR$, Let $P_i(c)$ be $\{(x_1,\dots,x_n)\in\RR^n\,:\,x_i=c\}$.
  Then, for all $c\in[0,z]$, $f\bigr|_{P_i(c)}$
  is bounded, otherwise, since $P_i(c)$ is homeomorphic
  to $\RR^{n-1}$, by induction $f\bigr|_{\Delta_{N\backslash\{i\}}(c)}$ is unbounded, and then
  by Lemmas \ref{lemma0} and \ref{lemma cof 1} $f$ is $N$-cofinal.
  Let $d(c)$ be a bound for $f\bigr|_{P_i(c)}$, $\{c_m\}_{m\in\omega}$ be
  a dense subset of $[0,z]$ and let $d=\sup_{m\in\omega}d(c_m)$. By density and continuity,
  $f\bigr|_{P_i(c)}$ is bounded by $d$ for all $c\in[0,z]$, and therefore
  $f\bigr|_{Q_i}$ is also bounded by $d$, contradiction.
\endproof

\begin{cor}\label{cor cof}
  The cofinality classes of continuous functions $\RR^n\to\RR$ are in bijection 
  with the antichains
  of $\mathcal{P}(N)\backslash\{\emptyset\}$.
\end{cor}

\proof 
  By Lemmas \ref{lemma cof 1}--\ref{lemma cof 2},
  the cofinality classes of continuous functions $\RR^n\to\RR$ are in
  bijection with the subsets $\mathfrak{I}$
  of $\mathcal{P}(N)\backslash\{\emptyset\}$
  satisfying the condition that 
  if $I\in\mathfrak{I}$ is non-empty and $J\supset I$, then $J\in\mathfrak{I}$.
  (Given such a $\mathfrak{I}$ it is easy to find $f$ such that $\mathfrak{C}(f)=\mathfrak{I}$,
  see below.)
  To any such $\mathfrak{I}$ corresponds bijectively
  an antichain given by its minimal
  elements. The empty antichain corresponds to bounded maps.
\endproof

Theorem \ref{thm2} follows immediately from this corollary and Theorem \ref{thm1}.
It is easy to find a representant for each class of maps $\RR^n\to\RR$:

\begin{defi}\label{canonical representants}
  Let $\mathscr{I}=\{I_1,\dots,I_k\}$ be an antichain in 
  $\mathcal{P}(N)\backslash\{\emptyset\}$. The canonical representant of the
  cofinality class $\mathscr{I}$ is given by
  \begin{equation}
    \label{eq:canonical representants}
    f_{\mathscr{I}}(x_1,\dots,x_n)\definition
    \left\{
    \begin{array}{l} \displaystyle
      \max_{\ell=1,\dots,k}\left\{\min_{i\in I_\ell} \{x_i\}\right\}
      \text{ if $\mathscr{I}\not=\emptyset$} \\
     \\
      \displaystyle
      0 \text{ if $\mathscr{I}=\emptyset$.}
    \end{array}
    \right.
   \end{equation}
\end{defi}
One checks easily that if $\mathscr{I}$ is an antichain,
$\mathscr{I}$ contains exactly the minimal elements of $\mathfrak{C}(f_{\mathscr{I}})$.
Assuming Theorems \ref{thm1}-\ref{thm2}, it is easy to show
that $[\RR^n,\RR^n]$ is isomorphic to a monoid of $(2^n-1)\times (2^n-1)$ matrices
with entries $0,1$.
\begin{defi}\label{direction matrix}
  Let $f:\RR^n\to\RR^n$ be continuous. We define its direction matrix
  $D(f)=(D_{I,J}(f))_{I,J\in \mathcal{P}(N)\backslash\{0\}}$ by
  $D_{I,J}(f)=1$ if there is some $c\in\RR^{n-|J|}$ such that $f(\Delta_I)\cap\Delta_J(c)$
  is unbounded in $\Delta_J(c)$, and $D_{I,J}(f)=0$ otherwise.
\end{defi}
If $D(f)=D(g)$, then $\mathfrak{C}(\pi_i\circ f)=\mathfrak{C}(\pi_i\circ g)$
for each $i\in N$, thus by Theorem \ref{thm1}, $f$ and $g$ are homotopic, the converse being obviously
true. Notice that by continuity, for a fixed $I$ there is at most one $J$ such that $D_{I,J}(f)=1$.
\begin{prop}\label{lemme monoide}
  If $f,g:\RR^n\to\RR^n$ are continuous,
  $D(f\circ g)=D(g)\cdot D(f)$.
\end{prop}
\proof
  By taking canonical representants for each coordinate of $f,g$,
  (Definition \ref{canonical representants}), we may assume that
  $f(\Delta_I)\cap\Delta_J$ is unbounded if and only if $D_{I,J}(f)=1$ (that is, the `$c$' of
  Definition \ref{direction matrix} is $0$), and similarly for $g$.
  The proof is then a routine check.
\endproof

Of course, not every $(2^n-1)\times (2^n-1)$ matrix of $0,1$ is a direction matrix, there are some
restrictions (which seem however quite tedious to describe).


\section{Partition properties}

The goal of this section is to prove an analog of Lemma 2.2 in \cite{meszigues+Dave} which says that
if $f:\RR\to\RR$ is continuous and cofinal, there is a partition
$\{P(\beta)\}_{\beta\in\omega_1}$ of $\RR$, with $P(\beta)=[x_\beta,x_{\beta+1}]$, such that
$f(P(\beta))=P(\beta)$ for each $\beta$. The homotopy question is then reduced to the
trivial problem of finding homotopies (here, between $f$ and the identity map)
defined in $P(\beta)$ and leaving $\partial P(\beta)=\{x_\beta,x_{\beta+1}\}$ fixed.
If $f:\RR\to\RR$ is bounded, then by Lemma \ref{lemmaintro2}, $f$ is constant outside a bounded
set, and the homotopy question is again trivial.

For maps $f:\RR^n\to\RR$ with $\mathfrak{C}(f)=\mathcal{P}(N)\backslash\{\emptyset\}$,
\cite[Lemma 2.2]{meszigues+Dave} may be applied to define homotopies in
$P(\beta)=[x_\beta,x_{\beta+1}]^n\backslash[0,x_\beta[^n$. Since
$\cup_{\alpha\in\omega_1}P(\beta)=\RR^n$, it suffices to glue together the homotopies to
obtain Theorem \ref{thm1} in this special case.
If $f$ is bounded then $f$ is trivially homotopic to the constant map $0$.
The problem is more difficult if $f$ is cofinal in some but not all $I\subset N$, but the idea
is always to find a partition $\{[x_\beta,x_{\beta+1}]\}_{\beta\in\omega_1}$
of $\RR$,
such that if $x\in\RR^n$ is ``between $x_\beta$ and $x_{\beta+1}$ in a cofinal direction $I$''
(see Definition \ref{A+A-} and (\ref{P}) below), then $f(x)\in [x_\beta,x_{\beta+1}]$.
Moreover, $f$ will be constant ``along bounded directions'' for $x$ sufficiently large (in these
directions), see Lemmas \ref{partition bounded}--\ref{BI and UI}.

\begin{lemma}\label{partition bounded}
   Let $f:\RR^n\to\RR$ be bounded. Then, there are $x,c\in\RR$ such that
   $f\bigr|_{[x,\omega_1[^n}\equiv c$.
\end{lemma}
\proof
   For $n=1$, this is Lemma \ref{lemmaintro2}. If $n\ge 1$, since $f\bigr|_{\Delta_N}$ is
   bounded, there are $(x,\dots,x)\in \Delta_N$ and $c\in\RR$ such that
   for all $x'\ge x$, $f(x',\dots,x')=c$.
   Let $c_m<c<c_m'$ be sequences converging to $c$. By Lemma \ref{structure open sets},
   since $f^{-1}(\,]c_m,c_m'[\,)\supset(\Delta_N\backslash[0,x[^n)$,
   there is $x_m$ such that $f^{-1}(\, ]c_m,c_m'[\, )\supset [x_m,\omega_1[^n$.
   Thus, $f^{-1}(c)=\bigcap_m f^{-1}(\, ]c_m,c_m'[\, )\supset [\sup_m x_m, \omega_1[^n$.
\endproof

It is useful to introduce the following notation:

\begin{defi}\label{MI}
  If $I\subset N$, $c\in\RR$, we set
  $M_I(c)\definition\{x\in\RR^n\,:\,x_I\in[c,\omega_1[^{|I|}\}$.
\end{defi}
\begin{lemma}\label{Lemma5.2}
 Let $I\subset N$ be non-empty.
 \begin{itemize}
 \item[1)]
  Suppose that $f$ is $I$-bounded.
  Then, for all $c\in\RR$, there is $d(c)$ minimal
  such that for all $b\in [0,c]^{n-|I|}$,
  $f$ restricted to
  \begin{equation}\label{def E}
    E_I(b,d(c))\definition\{x\in M_I(d(c))\,:\,x_{N\backslash I}=b\}
  \end{equation}
  is constant.
 \item[2)]
  Suppose that $f$ is $I$-cofinal. Then, for all $c\in\RR$, there is $\wt{d}(c)\in\RR$
  minimal such that for all $b\in[0,c]^{n-|I|}$,
  $f(\Delta_I(b)\cap M_I({\wt{d}(c)}))\subset[c,\omega_1[$.
 \end{itemize}
\end{lemma}
\proof
1)
  By Lemmas \ref{lemma0} and \ref{partition bounded}, for all
  $b\in [0,c]^{n-|I|}$ there is $d'(b)$ such that
  $f$ restricted to $E_I(b,d'(b))$ is constant.
  Let $\{b_m\}_{m\in\omega}$ be a dense subset
  of $[0,c]^{n-|I|}$. Then $d(c)=\sup_{m\in\omega}d'(b_m)$ has the required property.
  If $d(c)$ is not minimal, take the minimal one (which exists by continuity of $f$).

2)
  Let $b\in[0,c]^{n-|I|}$; since $\Delta_I(b)$ is homeomorphic to $\RR$,
  by \cite[Lemma 2.2]{meszigues+Dave} there is $\wt{d}'(b)$ such that
  $f(\Delta_I\cap M_I({\wt{d}'(b)}))\subset [c,\omega_1[$.
  As in 1), $\wt{d}(c)=\sup_{m\in\omega}(\wt{d}'(b_m))$ has the required property if
  $\{b_m\}_{m\in\omega}$ is dense in $[0,c]^{n-|I|}$.
\endproof

\begin{lemma}\label{BI and UI}
  Assume that $f$ is $I$-bounded.
  In the notations of Lemma \ref{Lemma5.2}, let $p_I^{\text{\rm bd}}(c)=\max\{c,d(c)\}$.
  Then,
  $p_I^{\text{\rm bd}}\bigr|_{\omega_1}$ is monotone increasing and continuous, and
  $\{c\in\RR\,:\,p_I^{\text{\rm bd}}(c)=c\}$ contains a closed unbounded set.
  If $f$ is $I$-cofinal, then $p_I^{\text{\rm cf}}(c)=\max\{c,\wt{d}(c)\}$ has the same properties.
\end{lemma}
It is not true in general that $p_I^{\text{cf}}$ and $p_I^{\text{bd}}$ are continuous in $\RR$. We will use
the fact that in $\omega_1$ we only have limits ``from below''.
\proof
Assume that $f$ is $I$-bounded.
We first prove that $p_I^{\text{bd}}$ is continuous (monotonicity is clear by definition).
It is enough to prove it for $d(c)$.
Let $\alpha_m\in\omega_1\subset\RR$ ($m\in\omega$)
be a sequence converging to $\alpha\in\omega_1$. We may assume that for each $m$,
$\alpha_m\le\alpha$. Let $d'(\alpha)=\lim_{m\to\infty}d(\alpha_m)$. By monotonicity of $d(\cdot)$,
the limit exists and $d'(\alpha)\le d(\alpha)$.
By minimality of $d(\alpha)$, it is enough to show that for each $b\in[0,\alpha]^{n-|I|}$,
$f$ restricted to $E_I(b,d'(\alpha))$ is constant.
\\
Assume for simplicity that $I=\{1,\dots,s\}$. Let $x=(x_1,\dots,x_n)\in E_I(b,d'(\alpha))$,
that is, $x_1,\dots,x_s\ge d'(\alpha)$ and $(x_{s+1},\dots,x_n)=b$.
For each $m\in\omega$, choose $b_m\in[0,\alpha_m]^{n-s}$ such that $b_m\to b$.
Put $x_m=(x_1,\dots,x_s,b_m)$. Since $d'(\alpha)\ge d(\alpha_m)$,
$x_m\in E_I(b_m,d(\alpha_m))$, so $f(x_m)=f(y_m)$ where
$y_m=(d(\alpha_m),\dots,d(\alpha_m),b_m)$.
By continuity, $f(x)=f(d'(\alpha),\dots,d'(\alpha),b)$. Since this holds for all $x$ in
$E_I(b,d'(\alpha))$, $f$ is constant on this set.
\\
We now prove that $K=\{ \alpha\in\omega_1\subset\RR\,:\,p_I^{\text{bd}}(\alpha)=\alpha\}$ is closed unbounded.
Closeness is immediate by continuity. Let $\alpha_0\in\omega_1$. Define inductively
$\alpha_m\in\omega_1$ ($m\in\omega$) such that $\alpha_m\ge p_I^{\text{bd}}(\alpha_{m-1})$.
By continuity, $\lim_{m\to\infty}\alpha_m=\alpha=p_I^{\text{bd}}(\alpha)\ge\alpha_0$. This
shows that $K$ is unbounded. The proof for $p_I^{\text{cf}}$ is similar.
\endproof

We have now taken care of the bounded directions. We proceed with the investigation of cofinal
directions. Lemma \ref{BI and UI} for $p_I^{\text{cf}}$ will be helpful.

\begin{defi}\label{A+A-}
  Let $I\subset{N}$ and $\alpha\in \omega_1$. We set:
  \begin{align*}
    A_I^-(\alpha)&\definition\left\{x\in\RR^n\,:\,
    \begin{array}{l}
      x_I\in[0,\alpha]^{|I|}\text{, and}  \\
      x_j\le |x_I|\,\forall j\in N
    \end{array}
    \right\}, 
    \\
    A_I^+(\alpha)&\definition\left\{x\in\RR^n\,:\,
    \begin{array}{l}
      x_I\in[\alpha,\omega_1[^{|I|}\text{, and}\\
      x_j\le |x_I|\,\forall j\in N
    \end{array}
    \right\}.
  \end{align*}
\end{defi}

The following lemma is the key argument for proving Theorem \ref{thm1}.

\begin{lemma}\label{partition lemma}
  If $f$ is $I$-cofinal, then
  $$
    \{\alpha\in \omega_1\,:\, f(A_I^+(\alpha))\subset[\alpha,\omega_1[\text{ and }
    f(A_I^-(\alpha))\subset[0,\alpha]\}
  $$
  is closed and unbounded.
\end{lemma}   
   
\proof
   1)
      We first show that
      $F^-=\{\alpha\in \omega_1\,:\,f(A_I^-(\alpha))\subset[0,\alpha]\}$ is closed
      and unbounded.
      Let $\alpha_m\to\alpha$, $m\in\omega$. One may assume $\alpha_m\le\alpha$.
      Since $\overline{\cup_{m\in\omega}A_I^-(\alpha_m)}=A_I^-(\alpha)$,
      by continuity $f(A_I^-(\alpha))\subset \overline{(f(\cup_{m\in\omega}A_I^-(\alpha_m))}
      \subset [0,\alpha]$. Thus, $F^-$ is closed.
      We now show that $\{\alpha\,:\,f([0,\alpha]^n)\subset[0,\alpha]\}$ is unbounded, which
      implies that $F^-$ is unbounded since $A_I^-(\alpha)\subset [0,\alpha]^n$. So, let 
      $\beta_0\in \omega_1$. For $m\in\omega$, we define $\beta_m\ge\beta_{m-1}$
      such that
      $f([0,\beta_{m-1}]^n)\subset[0,\beta_m]$ ($f$ being continuous, $f([0,\beta_{m-1}]^n)$
      is compact and thus bounded). Then, $\lim_{m\to\infty}\beta_m=\beta\in F^-$.

   2) We now show that
      $F^+=\{\alpha\in \omega_1\,:\,f(A_I^+(\alpha))\subset[\alpha,\omega_1[\}$ is closed
      and unbounded.
      The proof that $F^+$ is closed is like in 1), using
      $\cap_{m\in\omega}A_I^+(\alpha_m)=A_I^+(\alpha)$. To prove that $F^+$ is unbounded, we
      use the following claim:
      \begin{claim}
        For all  $\alpha\in \omega_1$, there is $\beta(\alpha)\ge\alpha$ such that
        $f(A_I^+(\beta(\alpha)))\subset[\alpha,\omega_1[$.
      \end{claim}
      This suffices to finish the proof: given $\alpha_0$, we define the sequence
      $\alpha_m=\beta(\alpha_{m-1})$ ($m\in\omega$) whose limit $\alpha\ge\alpha_0$ is in $F^+$.
      \proof[Proof of the claim.]
        To simplify, assume that $I=\{1,\dots,s\}$.
        Ab absurdo, suppose that:
        \begin{equation}\label{ab absurdo}
          \forall\alpha,\forall\beta\ge\alpha,\,\exists
          x_\beta\in A_I^+(\beta)\text{ with }f(x_\beta)\le\alpha.
        \end{equation}
	For each $\beta$, we fix $x_\beta\in A_I^+(\beta)$
	such that (\ref{ab absurdo}) holds. We now proceed
	in several steps.\\
        a) We first show that
        $$
         \Gamma=\left\{\gamma\,:\,
         \begin{array}{l}
           \exists b_{s+1},\dots,b_n\in\RR \text{ such that} \\
           f(\gamma,\dots,\gamma,b_{s+1},\dots,b_n)\in[0,\alpha]
         \end{array}
         \right\}
        $$
        is closed and unbounded. Indeed, given $\gamma_m\to\gamma$, $\gamma_m\in\Gamma$,
        there corresponds sequences $b_{j,m}$ $j=s+1,\dots,n$. Taking convergent
        subsequences, we see that $\gamma\in\Gamma$, which is thus closed.
        Now, given $\beta\in \omega_1$, we may define
        $\gamma_0=\max\{\alpha,\beta\}$. Then, by induction, choose
        $\gamma_m\ge |x_{\gamma_{m-1}}|$; thus
        $f(x_{\gamma_m})\le\alpha$ (recall that each $x_{\beta}$ satisfies (\ref{ab absurdo})).
        Taking a convergent subsequence of the $x_{\gamma_m}$, we obtain an
        $x=(\gamma,\dots,\gamma,b_{s+1},\dots,b_n)$ with $\gamma\ge\beta$
        and $f(x)\le\alpha$, showing that $\Gamma$ is unbounded.
        \esp
        b) We then set:
        $$
          C(\gamma)\definition\left\{b=(b_{s+1},\dots,b_n)\in\RR^{n-s}\,:\,
          \begin{array}{l}
            \exists \beta\ge\gamma\text{ with} \\
            f(\beta,\dots,\beta,b)\le\alpha
          \end{array}
          \right\}.
        $$
        $C(\gamma)$ is nonempty,
        closed, and if $\gamma'\ge\gamma$, $C(\gamma')\subset C(\gamma)$.
        \esp
        c) For all $\gamma\in\omega_1$, $C(\gamma)$ is unbounded. Indeed, Lemma \ref{BI and UI}
	(for $p_I^{\text{cf}}$) implies that the set of $\beta$ satisfying
        $f(\Delta_I(b)\cap M_I(\beta))\subset[\beta,\omega_1[$
	for all $b\in[0,\beta]^{n-s}$ is unbounded.
	Thus, if $C(\gamma)$ is bounded, $C(\gamma)\subset [0,\beta]^{n-s}$
	for one such $\beta >\alpha$, and then $C(\beta)\subset [0,\beta]^{n-s}$, which
	implies that $C(\beta)$ is empty, contradicting b).
	\esp
        d) By c) and Lemma \ref{structure closed sets}, for all $\gamma\in\omega_1$ the set
        $$
          \mathfrak{J}(\gamma)=\left\{J\subset \{s+1,\dots,N\}\,:\,
          \begin{array}{l}
             \exists c\in\RR^{n-s-|J|}\text{ such that } \\
             \Delta_{J}(c)\cap C(\gamma)
             \text{ is unbounded}
          \end{array}
          \right\}
        $$
        is nonempty. Since $\mathfrak{J}(\gamma')\subset\mathfrak{J}(\gamma)$ if
	$\gamma'>\gamma$, $\cap_{\gamma\in \omega_1}\mathfrak{J}(\gamma)$ is also
        nonempty, let $J$ be in this intersection. We may assume
        $J=\{s+1,\dots,s_1\}$. So, for all $\gamma\in\omega_1$, and all $x\in\RR$, there
        is $y\ge x$ and $c_{s_1+1},\dots,c_n$ such that
        $(y,\dots,y,c_{s_1+1},\dots,c_n)\in C(\gamma)$. In other words,
        for any $\gamma$, there are $y\ge\gamma$, $\beta\ge\gamma$ and
        $c_{s_1+1},\dots,c_n$ such that
        \begin{equation}\label{beta et y}
          f(\underbrace{\beta,\dots,\beta}_{s},\underbrace{y,\dots,y}_{s_1-s},
          c_{s_1+1},\dots,c_n)\le\alpha.
        \end{equation}
        Given $\gamma_0$, we may define sequences $\gamma_m,\beta_m,y_m,c_{s_1+1,m},\dots,c_{n,m}$
	by letting
        $\gamma_m\ge\max\{\gamma_{m-1},y_{m-1}\}$ and choosing
        $y_m\ge \gamma_m$, $\beta_m\ge\gamma_m$ and $c_{s_1+1,m},\dots,c_{n,m}$
	satisfying (\ref{beta et y}). We thus have
        $\gamma_m\ge\beta_{m-1}\ge\gamma_{m-1}$ and $\gamma_m\ge y_{m-1}\ge\gamma_{m-1}$,
        these three sequences converge to the same point $\gamma$. Taking
        convergent subsequences of $c_{s_1+1,m},\dots,c_{n,m}$,
        we have found an $x=(\gamma,\dots,\gamma,c_{s_1+1},\dots,c_n)\in\RR^n$ such that
        $\gamma\ge\gamma_0$ and $f(x)\le\alpha$.
        That is, we proved that
        $$
         \Gamma_1=\left\{\gamma\,:\,
         \begin{array}{l}
           \exists b_{s_1+1},\dots,b_n\in\RR \text{ such that}\\
           f(\gamma,\dots,\gamma,b_{s_1+1},\dots,b_n)\in[0,\alpha]
         \end{array}
         \right\}
        $$
        is unbounded, and its closeness is immediate.
        \esp
        e) We may thus go back to b) with $s_1$ instead of $s$, and proceed by induction
	until we obtain
        $\{\gamma\,:\,f(\gamma,\dots,\gamma)\le\alpha\}$ is closed and unbounded.
	(In c), we use Lemma \ref{lemma cof 1} which ensures that $f$ is $\{1,\dots,s_1\}$-cofinal.)
        But by Lemma \ref{BI and UI} for $p_N^{\text{\em cf}}$,
	this implies that $f\bigr|_{\Delta_N}$ is bounded,
        which is contradicts Lemma \ref{lemma cof 1} since $f$ is
        $I$-cofinal. This proves the claim.

\endproof


\section{Proof of the main theorem}

We will now use our partition properties to define an homotopy between $f$ and the canonical
representant of its cofinality class.
We first recall the following triviality:
\begin{lemma}\label{trivial homotopy lemma}
  Let $g,h:X\to Y$ be continuous and $Y$ be homeomorphic to $[0,1]^d$. Then,
  there is a homotopy $\phi_t$ such that $\phi_0=g,\phi_1=h$ and for all $t$,
  $\phi_t\bigr|_Q=id$, where $Q=\{x\in X\,:\,f(x)=g(x)\}$.
\end{lemma}
\proof
  Let $\varphi:Y\to[0,1]^d$ be an homeomorphism. Then,
  $\phi_t(x)=\varphi^{-1}\Bigl(\varphi(g(x))\cdot(1-t)+\varphi(h(x))\cdot t\Bigr)$ has the required properties.
\endproof

\begin{figure}[h]
    \begin{center}
    \epsfig{figure=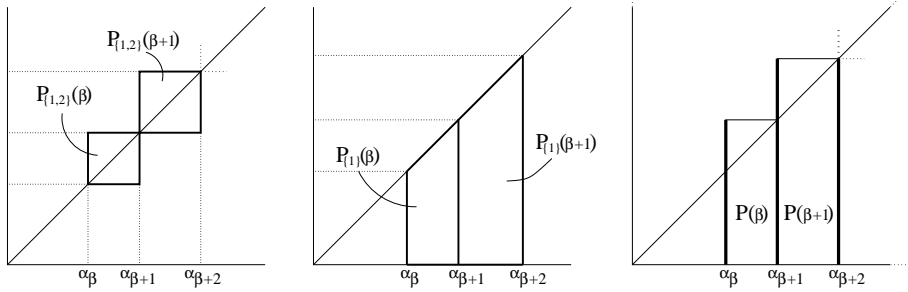,  height=4cm}
    \caption{The sets $P_I(\alpha)$ ($I={\{1\}},{\{1,2\}}$; $\alpha=\beta,\beta+1$) for $h$.}
    \label{fig:1}
    \end{center}
  \end{figure}

\begin{figure}[h]
    \begin{center}
    \epsfig{figure=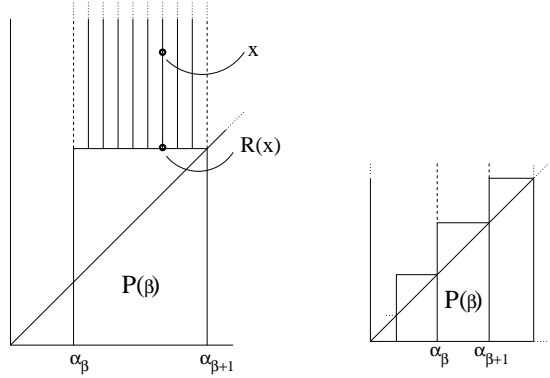,  height=5cm}
    \caption{Constant verticals and $R(x)$ for $h$.}
    \label{fig:2}
    \end{center}
  \end{figure}

\proof[Proof of Theorem 1]
  Let $f:\RR^n\to\RR$ be continuous.
  The case $\mathfrak{C}(f)=\emptyset$ (and thus $f$ bounded) being trivial, we may assume
  $\mathfrak{C}(f)\not=\emptyset$.
  To clarify the exposition, we fix some map $h:\RR^2\to\RR$ with
  $\mathfrak{C}(h)=\{\{1\},\{1,2\}\}$ to serve us as an example. We shall carry the proof
  for general $f:\RR^n\to\RR$ and for $h$ together.
  Let $\mathscr{J}$ be the minimal elements of $\mathfrak{C}(f)$.
  We shall show that $f$ and $f_{\mathscr{J}}$ (defined by (\ref{eq:canonical representants}))
  are homotopic. In the case of $h$, $\mathscr{J}=\{\{1\}\}$ and $f_{\mathscr{J}}$ is the projection
  on the first coordinate.
  \\
  By Lemmas \ref{lemmaintro1}, \ref{BI and UI} and \ref{partition lemma}, the set $\Theta$ of ordinals
  $\alpha$ satifying
  \begin{align}
    \label{part.1}
    \forall I\in\mathfrak{C}(f),\quad &f(A_I^+(\alpha))\subset[\alpha,\omega_1[,\,
    f(A_I^-(\alpha))\subset[0,\alpha],\,
    \\
    &\text{and}\nonumber\\
    \label{part.2}
   \forall J\notin\mathfrak{C}(f),\quad & \alpha=p_J^{\text{\rm bd}}(\alpha)
  \end{align}
  is closed and unbounded.
  For all $\beta\in \omega_1$, we then choose $\alpha_\beta\in\Theta$ as follows:
  \begin{align*}
    \alpha_{\beta+1}&=\min\Theta\cap [\alpha_\beta+1,\omega_1[, \\
    \alpha_\beta & =\sup_{\gamma<\beta}\alpha_\gamma\text{ if $\beta$ is a limit ordinal.}
  \end{align*}
  Then, for all $I\in\mathfrak{C}(f)$, we set $P_I(\beta)=A_I^+(\alpha_\beta)\cap
  A_I^-(\alpha_{\beta+1})$, that is,
  $P_I(\beta)=\{x\in\RR^n\,:\, x_I\in[\alpha_\beta,\alpha_{\beta+1}]^{|I|} \text{ and }
  x_j\le |x_I| \, \forall j\in N\}$. Finally, set
   \begin{equation}\label{P}
    P(\beta)=\bigcup_{I\in\mathfrak{C}(f)}P_I(\beta).
   \end{equation}
  The corresponding sets for $h$ are shown on Figure \ref{fig:1}.
  By (\ref{part.1}), for $I\in\mathfrak{C}(f)$ and all $\beta$,
  $f(P_I(\beta))=[\alpha_\beta,\alpha_{\beta+1}]$. Notice that we also have
  $f_{\mathscr{J}}(P_I(\beta))=[\alpha_\beta,\alpha_{\beta+1}]$.
  One easily checks that
  $P_I(\beta)\cap P_I(\beta+1)=\{x\in\RR^n\,:\, x_i=\alpha_{\beta+1}\,\text{for }i\in I\text{ and }
  x_j\le\alpha_{\beta+1}\,\forall j\in N\}$, and thus
  $f(P_I(\beta)\cap P_I(\beta+1))=\{\alpha_{\beta+1}\}$;
  if
  $I,J\in\mathfrak{C}(f)$, then $P_I(\beta)\cap P_J(\beta)\subset P_{I\cup J}(\beta)$
  (recall that by Lemma \ref{lemma cof 1} $I\cup J\in\mathfrak{C}(f)$), and
  $P_I(\beta)\cap P_J(\beta+1)\subset P_J(\beta)\cap P_J(\beta+1)$.
  Thus,
  \begin{align*}
    P(\beta)\cap P(\beta+1)&=\bigcup_{I\in\mathfrak{C}(f)} (P_I(\beta)\cap P_I(\beta+1))
    \\
    &= \bigcup_{I\in\mathfrak{C}(f)}
    \left\{x\in\RR^n\,:\,
    \begin{array}{l}
      x_i=\alpha_{\beta+1}\,\text{for}\,i\in I,\\
      \,x_j\le\alpha_{\beta+1}\,\forall j\in N
    \end{array}
    \right\},
    \end{align*}
  and
  \begin{equation}\label{P intersection}
     f(P(\beta)\cap P(\beta+1))=\{\alpha_{\beta+1}\}=f_{\mathscr{J}}(P(\beta)\cap P(\beta+1)).
  \end{equation}
  For $h$, this means that the bold vertical boundaries of the rightmost picture in
  Figure \ref{fig:1} are ``projected
  down''.
  \\
  We can now apply Lemma \ref{trivial homotopy lemma} with $X=P(\beta)$ and
  $Y=[\alpha_\beta,\alpha_{\beta+1}]$ to find homotopies $\phi^\beta_t$ defined on $P(\beta)$
  such that $\phi^\beta_0=f\bigr|_{P(\beta)}$ and $\phi^\beta_1=f_{\mathscr{J}}\bigr|_{P(\beta)}$.
  By (\ref{P intersection}), if $x\in P(\beta)\cap P(\beta+1)$,
  $\phi^\beta_t(x)=\phi^{\beta+1}_t(x)=\alpha_{\beta+1}$ for all $t\in[0,1]$, and we can ``glue''
  together the $\phi^\beta_t$ to obtain an homotopy $\phi_t$ between $f$ and $f_{\mathscr{J}}$
  defined on $P=\bigcup_{\beta\in\omega_1}P(\beta)\subset\RR^n$. We shall now explain how to
  extend this homotopy to all of $\RR^n$.\\
  First, consider our example $h$.
  We have depicted the situation of $P(\beta)$ in Figure \ref{fig:2}.
  By Lemma (\ref{part.2}),
  $h$ restricted to any vertical
  line depicted in Figure \ref{fig:2} is constant
  (these vertical lines are exactly the $E_{\{2\}}(b,\alpha_{\beta+1})$
  for $b\in [\alpha_\beta,\alpha_{\beta+1}]$). We can then define $R(x)\in\partial P$
  for $x\notin P$
  as in this figure, and then $h(x)=h(R(x))$.
  Since the vertical boundaries of $P(\beta)$ are both mapped by $h$ to one point,
  the ambiguity of the definition of $R(x)$ for
  $x$ lying on one of the dashed lines of
  Figure \ref{fig:2} does not cause any trouble.
  If we extend $R$ by the identity in $P$, $R(x)$ will be non-continuous (due to the
  above ambiguities), but
  $\wt{\phi}_t(x)=\phi_t((R(x))$ is continuous, and is then an homotopy between
  $h$ and $f_{\{1\}}$ (since $f_{\{1\}}(x)=f_{\{1\}}(R(x))$ as well).
  One sees easily that for fixed $t$, $\wt{\phi}_t$ it is constant on the verticals depicted
  and fixes $P(\beta)\cap P(\beta+1)$ for all $t$.
  Moreover, limit ordinals $\beta$ do not cause any trouble.
  (Strictly speaking, the homotopy is not defined in $[0,\alpha_0]\times\RR$, but we may squish
  this set continuously to $\{\alpha_0\}\times\RR$ and no bother about it.)

  Let us do the general case, i.e. go back to $f:\RR^n\to\RR$.
  The only difference is a heavier formalism and no new idea is needed, we shall thus pass
  quickly over the details.
  Let $x\in\RR^n$, with coordinates
  $x_{i_1}\ge x_{i_2}\ge\cdots\ge x_{i_n}$. Choose $\beta$ such that
  $x_{i_1}\in[\alpha_\beta,\alpha_{\beta+1}]$. (The fact that $\beta$ is not always unique
  is not important.) Let $k$ be maximal such that
  $x_{i_1},x_{i_2}\dots,x_{i_k}\in [\alpha_\beta,\alpha_{\beta+1}]$. If
  $I=\{i_1,\dots,i_k\}\in\mathfrak{C}(f)$, then $x\in P_I(\beta)\subset P$.
  If $I\notin\mathfrak{C}(f)$, choose $\ell\ge k$ maximal such that
  $I'=\{i_1,\dots,i_\ell\}\notin\mathfrak{C}(f)$ and for some $\beta'$,
  $x_{i_1},\dots,x_{i_\ell}\ge\alpha_{\beta'}$. Then, by definition,
  $x\in E_{I'}(b,\alpha_{\beta'})$ where
  $b=(x_{i_{\ell+1}},\dots,x_{i_n})\in[0,\alpha_{\beta'}]^{n-\ell}$ (see (\ref{def E})).
  We can choose $\beta'$ minimal, and then $\beta'$ is successor if $x_{i_{\ell+1}}<\alpha_{\beta'}$.
  Then,
  by minimality of $\beta'$, $x_{i_{\ell+1}}\in\,]\alpha_{\beta'-1},\alpha_{\beta'}[\,$.
  Set then $R(x)=(R_1(x),\dots,R_n(x))$ by letting $R_{i_s}(x)=x_{i_s}$ for $s=1,\dots,\ell$ and
  $R_{i_s}(x)=\alpha_{\beta'-1}$ for $s=\ell+1,\dots,n$.
  By maximality of $\ell$, $J=\{i_1,\dots,i_{\ell+1}\}\in\mathfrak{C}(f)$, and
  $R(x)\in\partial P_J(\beta'-1)$. If
  $x_{i_{\ell+1}}=\alpha_{\beta'}$, we set
  $R_{i_s}(x)=\alpha_{\beta'}$ for $s=\ell+1,\dots,n$, and then $R(x)\in\partial P_J(\beta')$.
  By (\ref{part.2}) and Lemma \ref{BI and UI} for $p_I^{\text{bd}}$, $f(x)=f(R(x))$ in both cases.
  Extending $R(x)$ to all $\RR^n$ by $R(x)=x$ for $x$ in $P$, we may define
  $\wt{\phi}_t(x)=\phi_t(R(x))$ and check as in $h$'s case that $\wt{\phi}_t$ is continuous
  and sends $f$ to $f_{\mathscr{J}}$.
\endproof


\section{Other manifolds}
In this section we consider some other non-metrizable manifolds and state some theorems about their homotopy classes.

Recall first that the long line $\LL$ is the union of two copies $\LL^+,\LL^-$ of $\RR$ glued
at $0$. In order to code maps $\LL^n\to\RR$, we
let $N^\pm$ be the set of ``signed'' coordinates $\{+1,\dots,+n,-1,\dots,-n\}$
and say that $I\subset N^\pm$ is an admissible subset of $N^\pm$, which we denote by
$I\adm N^\pm$, if for all $i\in N$, $I$ does not contain both $+i$ and $-i$.
We then set $\mathcal{P}^a(N^\pm)=\{I\adm N^\pm\}$.
We have the following result:

\begin{thm}\label{thm3}
  $[\LL^n,\RR]$ is in bijection with the antichains of
  $\mathcal{P}^a(N^\pm)\backslash\{\emptyset\}$, and
  $[\RR^n,\LL]$ is the union of $[\RR^n,\LL^+]$ and $[\RR^n,\LL^-]$ where bounded maps
  in $\LL^+$ and $\LL^-$ are identified.
\end{thm}
\proof
The assertion about $[\LL^n,\RR]$ is proved as Theorem \ref{thm2}. For $[\RR^n,\LL]$, notice
that a continuous map $\RR\to\LL$ cannot be unbounded in both $\LL^+$ and $\LL^-$.
Thus, if $f\restrict{\Delta_N}$ is cofinal in $\LL^+$, $f\restrict{\Delta_I}$
cannot be cofinal in $\LL^-$ by Lemma \ref{lemma cof 1}, and the result
follows thus from Theorem \ref{thm2}.
\endproof

The homotopy classes of maps $\LL^n\to\LL$ can be classified as well, but
are harder to describe.

\begin{figure}[h]
    \begin{center}
    \epsfig{figure=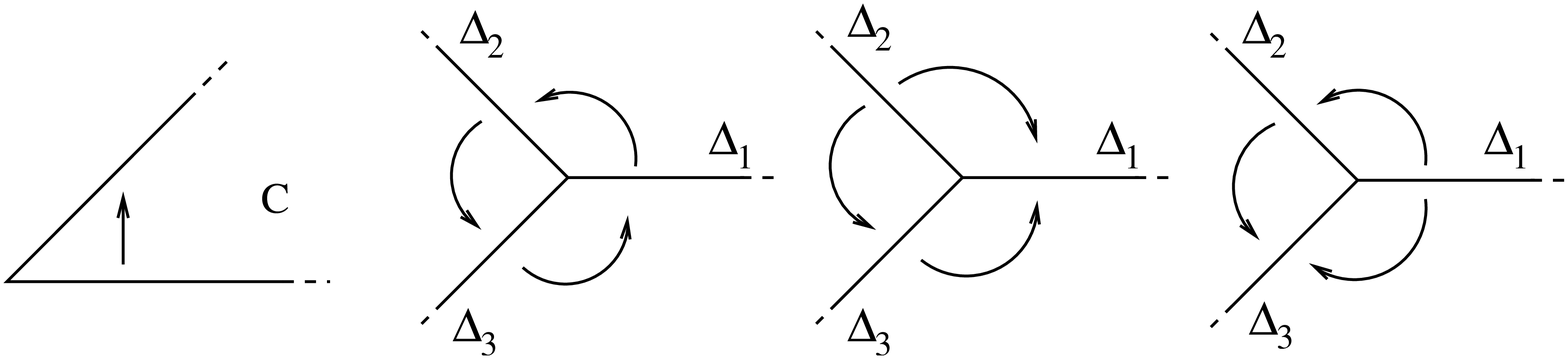,  width=10cm}
    \caption{$\mathsf{P}_S$ for $S=\{\uparrow\uparrow\uparrow\},\{\downarrow\uparrow\uparrow\},\{\uparrow\uparrow\downarrow\}$.}
    \label{fig:3}
    \end{center}
  \end{figure}

Let us give a few more examples in dimension $2$. Let $C$ be the set $\{(x,y)\in\RR^2:y\le x\}$.
Fix $k\in\omega$, and set $\wb{i}=i\mod (k+1)$ for $i\in\omega$.
Given a finite sequence $S=(s_1,\dots,s_k)$ of symbols $\uparrow$ and $\downarrow$, we define
the $S$-pipe
\begin{equation}\label{S pipe}
  \mathsf{P}_S=\bigcup_{i=1}^k C\times\{i\} / \sim
\end{equation}
where $x\sim y$ iff $x=y$ or
\begin{align*}
   x=((u,u),i),\, y=((u,0),\wb{i+1}) \text{ and } &s_i=s_{\wb{i+1}}=\uparrow ,&\text{ or}\\
   x=((u,u),i),\, y=((u,u),\wb{i+1}) \text{ and } &s_i=\uparrow,\, s_{\wb{i+1}}=\downarrow ,&\text{ or}\\
   x=((u,0),i),\, y=((u,0),\wb{i+1}) \text{ and } &s_i=\downarrow,\, s_{\wb{i+1}}=\uparrow ,&\text{ or}\\
   x=((u,0),i),\, y=((u,u),\wb{i+1}) \text{ and } &s_i=s_{\wb{i+1}}=\downarrow ,\\
\end{align*}
For instance,
$\LL^2=\mathsf{P}_{(\uparrow\downarrow\uparrow\downarrow\uparrow\downarrow\uparrow\downarrow)}$.
See Figure \ref{fig:3} for other examples.
Any such $\mathsf{P}_S$ is a long pipe if we remove one point (see \cite[Def 5.2]{Nyikos:1984}).
It is not difficult to
see that $\mathsf{P}_S$ and $\mathsf{P}_{S'}$
are homeomorphic iff $S$ and $S'$ differ only by a circular permutation
and/or a uniform exchange of $\uparrow$ and $\downarrow$. We may also prove that there
are exactly $k$ homotopy classes of unbounded maps $\RR\to \mathsf{P}_S$, canonical representants
being given by $x\mapsto \pi((x,0),i)$ if $s_i=\uparrow$ and $x\mapsto \pi((x,x),i)$
if $s_i=\downarrow$ ($\pi$ denotes the projection
of $\bigcup_{i=1}^k C\times\{i\}$ on the quotient space $\mathsf{P}_S$). We denote (the image of) these
representants by $\Delta_1,\dots,\Delta_k$ (see Figure \ref{fig:3}), and define $i$-cofinality and
$i$-boundedness of maps $\mathsf{P}_S\to\RR$ as in Section 2. $S$ gives a partial order $\prec$ on
$\{1,\dots,k\}$ defined by $i\prec j$ iff
$j=\wb{i+\ell}$ for some $\ell\le k$ and
$s_{\wb{i}} = \cdots = s_{\wb{i+\ell-1}}=\uparrow$, or
$i=\wb{j+\ell}$ for some $\ell\le k$ and
$s_{\wb{j+1}}= \cdots = s_{\wb{j+\ell}}=\downarrow$ (in other words, if we can pass from $i$ to $j$ following arrows).
It is possible that $i\prec i$, if all $s_j$ are equal.
As in
Lemmas \ref{lemma cof 1}--\ref{lemma cof 2}, we can prove
that if $f$ is $i$-cofinal and $i\prec j$, then $f$ is $j$-cofinal, and that
an unbounded map is $j$-cofinal for some (maximal) $j$. Applying the technique of
Theorems \ref{thm1}--\ref{thm2}, we obtain the following result:
\begin{thm}\label{thm4} 
  With the above notations, 
  $[\mathsf{P}_S,\RR]$ is in bijection with the antichains of the partially ordered set $\bigl<\{1,\dots,k\},\prec\bigr>$.
\end{thm}

For general non-metrizable manifolds, a complete description of its homotopy classes may be
quite difficult, even if all $\Pi_i(M)$ are zero, since
it may happen that  $[M,M]$, or even
$[M,\RR]$, is infinite (for instance, glue an infinite number of $C$ together in the same
fashion as the $\mathsf{P}_S$ above).
In all generality, what we can say is for instance the following.

\begin{prop}\label{contractible}
  Let $M$ be a manifold such that $\omega_1\subset M$
  (that is, there is an embedding $e:\omega_1\to M$). Then, $M$ is not contractible.
\end{prop}
It may be interesting to see if we can
weaken the hypotheses to ``$M$ is non-metrizable'', since there
are many non-metrizable manifolds that do not contain $\omega_1$. For instance, there are
smoothings of $\RR$ such that the tangent bundle with the $0$ section
removed does not contain any copy of $\omega_1$ (see \cite[class 7, p. 158]{Nyikos:1992}).
Notice that the assumption that $M$ is a manifold is essential: the cone over $\omega_1$ is
contractible.
\proof
We do not make the distinction between $\alpha\in\omega_1$ and $e(\alpha)\in M$, denoting
both by $\alpha$, and identify
$\omega_1$ and $e(\omega_1)$.
Suppose that there is a continuous $h:[0,1]\times X\to X$ with $h_0=id$ and $h_1\equiv y$ for
some $y\in M$.
Let $U$ be a chart around $y$. Notice that $\omega_1\not\subset U$.
Then, there is $s<1$ such that for all $s<t\le 1$,
$h_t(\omega_1)\subset U$.
(Otherwise take a sequence $t_m\to 1$, $m\in\omega$.
For each $m$ there is $\alpha_m\in\omega_1$ with
$h_{t_m}(\alpha_m)\notin U$. By taking a convergent subsequence of the
$\alpha_m$, we obtain $\alpha\in\omega_1$ for which
$h_1(\alpha)\notin U$, contradiction.) \\
Then, for $t>s$, $h_t|_{\omega_1}$ is eventually constant, i.e. $\exists\alpha\in\omega_1$
such that
$\forall \beta\ge\alpha$, $h_t(\beta) = h_t(\alpha)$
($U$ being homeomorphic to some bounded open set in $\R^d$, we may apply a modified version of
Lemma \ref{lemmaintro2}.)
We say that $h_t|_{\omega_1}$ is $\alpha$-eventually constant.
Let:
\begin{equation*}
  \tau\definition\inf\{t\in[0,1]\,\,:\,\, h_t|_{\omega_1}\text{ is eventually constant}\}.
\end{equation*}
We saw that $\tau<1$. There are two possibilities.

\begin{itemize}

\item[1)]
There is $\alpha$ such that $h_\tau|_{\omega_1}$ is
$\alpha$-eventually constant. Then, since $h_0=id$, $\tau>0$.
Choose a sequence $t_m\nearrow\tau$, $m\in\omega$. Since
$h_{t_m}|_{\omega_1}$ is non-eventually constant, if $V$ is a chart around $h_\tau(\alpha)$,
$\exists\beta_n\ge\alpha$ with $h_{t_n}(\beta_n)\notin V$. 
(Otherwise, like before,  $h_{t_n}|_{\omega_1}$ would be eventually constant.)
Taking a convergent subsequence of the $\beta_m$, 
we obtain $\beta\ge\alpha$ with $h_\tau(\beta)\notin V$, contradiction with $h_\tau(\beta)=h_\tau(\alpha)$.

\item[2)]
$h_\tau|_{\omega_1}$ is not eventually constant. Since $\tau<1$, let
$t_m\searrow\tau$, $m\in\omega$.
For each $m$, there is $\alpha_m$ with $h_{t_m}|_{\omega_1}$ $\alpha_m$-eventually constant.
Taking a subsequence converging to $\alpha$, we obtain that $h_\tau|_{\omega_1}$ is
$\alpha$-eventually constant, contradiction.

\end{itemize}
\esp
Therefore, such an $h_t$ cannot exist and $M$ is not contractible.
\endproof

{\it Acknowledgements.} I wish to thank David Cimasoni and Patrick Lapp.

\vskip1cm
\noindent
Mathieu Baillif \\
Section de Math\'ematiques\\
2-4 rue du Li\`evre\\
1211 Gen\`eve 24\\
Switzerland\\
{\tt Mathieu.Baillif@math.unige.ch}

\end{document}